\documentclass[11pt]{article}
\usepackage{amssymb}

\begin{document}

\def\R{{\mathbb R}}
\def\C{{\mathbb C}}

\title{On an example of a
transition from chaos to integrability
for magnetic geodesic flows}
\author{I.A. Taimanov
\thanks{Institute of Mathematics, 630090 Novosibirsk, Russia;
e-mail: taimanov@math.nsc.ru}}
\date{}
\maketitle

In this paper we give an example of a real-analytic Hamiltonian
system whose restrictions onto different energy levels are Anosov
flows for $E > E_{\mathrm{cr}}$, have only transitive trajectories
for $E=E_{\mathrm{cr}}$ and are analytically integrable for
$E<E_{\mathrm{cr}}$ where $E_{|mathrm{cr}}$ is some critical level
of energy. This gives a simple example of a transition from chaos
to integrability (on a fixed energy level) via a passing of the
energy through the critical level.

We have

\vskip3mm

{\bf Theorem.}\  \
{\sl Let $M$ be a closed oriented two-dimensional manifold with a metric of
constant curvature $K=-1$ and
let $d\mu$ be the volume form corresponding to this metric.
Then the motion of a charge particle on $M$ in a magnetic field given
by the form $d\mu$ is completely integrable in terms of
real-analytic integrals of motion on the energy levels $E < \frac{1}{2}$.}

\vskip3mm

{\sc Remark} 1. Hedlund proved that on the energy level
$E=\frac{1}{2}$ any trajectory of this flow, i.e., magnetic
geodesic flow, is transitive which means that it is everywhere
dense \cite{Hedlund}. On every energy level $E
> \frac{1}{2}$ the magnetic geodesic flow is conjugate to the
geodesic flow of a Finsler metric which is an Anosov flow
\cite{CIPP}.

{\sc Definition.} Given a Hamiltonian system on a $2n$-dimensional
symplectic manifold, we say that it is {\it integrable on the
fixed energy level} $\{H=E\}$ if there exist $(n-1)$ integrals of
motion $I_1,\dots,I_{n-1}$ such that they are defined in some
neighborhood of this level, where they commute with respect to the
Poisson brackets, and they are functionally independent on the
full measure subset of the energy level $\{H=E\}$ (see, for
instance, \cite{T,BT}).

Following the paper by Novikov \cite{Novikov},
recall that the motion of a charge particle on a Riemannian manifold $M$
in a magnetic field given by some
closed $2$-form $F$ is described by the Euler--Lagrange equations for
the Lagrangian
$$
L(x,\dot{x}) = \frac{1}{2} g_{ik}\dot{x}^i\dot{x}^k + A_i\dot{x}^i
$$
where $g_{ik}$ is a metric tensor and a $1$-form $A$ is locally defined as
$A = d^{-1}F$, i.e., trajectories are extremals of the action $\int L dt$.
If the form $F$ is exact then we say
the magnetic field is exact and in this case $A$ is globally defined.
Since the magnetic field comes into the
Euler--Lagrange equations via the form
$F$, these equations and their solutions (i.e., magnetic geodesics) do not
depend on a choice of $A = d^{-1}F$. In particular, this leads to an
introduction of multivalued functionals on spaces of closed curves for
looking for closed extremals by using analogs of the Morse theory
\cite{Novikov}.

A magnetic field does not contribute to the energy which equals
$$
E = \dot{x}^i \frac{\partial L}{\partial \dot{x}^i} - L
= \frac{1}{2}g_{ik}\dot{x}^i\dot{x}^k.
$$
Thus in the Hamiltonian formalism the magnetic geodesic flow is
described by Hamiltonian equations on the cotangent bundle
$T^\ast M$ with the same Hamiltonian function
$H(x,p) = \frac{1}{2}g^{ik}p_i p_k$ as the geodesic flow but with respect
to another symplectic structure. In this event it is given by a closed form
$\omega = dp_i \wedge dx^i + \pi^\ast F$ where $\pi: T^\ast M \to M$ is the
natural projection.

On a fixed energy level $E$ magnetic geodesics are
extremals of another action functional
$\int L_E dt$ where the Lagrangian $L_E$
takes the form
\begin{equation}
\label{1}
L_E(x,\dot{x}) = \sqrt{2E g_{ik} \dot{x}^i \dot{x}^k} + A_i \dot{x}^i =
\sqrt{2E} \left(\sqrt{g_{ik} \dot{x}^i \dot{x}^k} + \frac{1}{\sqrt{2E}}
A_i \dot{x}^i\right).
\end{equation}

Consider the magnetic geodesic flow corresponding to
a hyperbolic metric, i.e., with constant negative curvature $K=-1$,
and its volume form $d\mu$ on a compact surface.
It follows from (\ref{1}) that in this case
magnetic geodesics on the energy
level $E$ are lines of constant geodesic curvature
$$
k_g = \frac{1}{\sqrt{2E}}
$$
passed in the clockwise direction.
These lines are described by solutions to
the Euler--Lagrange equations for the multivalued
Lagrangian
$$
\sqrt{g_{ik} \dot{x}^i \dot{x}^k} + \frac{1}{\sqrt{2E}}
A^\alpha_i \dot{x}^i
$$
with $dA^\alpha = F$ in a domain $U_\alpha$ with
$[F\vert_{U_\alpha}] = 0 \in H^2(U_\alpha,\R)$ (we can take for $U_\alpha$
the complement to any point which does not lie on the line).

Hedlund described all lines with constant geodesic curvature on the
hyperbolic plane in \cite{Hedlund}. Their classification is as follows.
Let us realize the hyperbolic plane as the disc $D = \{|z| < 1\}$, on the
complex plane $\C$, with the metric
$$
ds^2 = \frac{dzd\bar{z}}{(1-|z|^2)^2}.
$$
The group of isometries is $SU(1,1)/\pm 1$ which acts as
$$
z \to \frac{az+b}{cz+d}, \ \ \
\left(\begin{array}{cc} a & b \\ c & d
\end{array}\right) \in SU(1,1).
$$
Then the lines of the constant geodesic curvature $k_g$ are

\begin{itemize}
\item
$k_g=0$: geodesics, which in the Euclidean metric on $\C$ are arcs of circles
orthogonal to $\partial D$;

\item
$0 < k_g < 1$: hypercycles, which are arcs of Euclidean circles meeting
$\partial D$ in two different points. All points of a hypercycle are
equidistant (in the hyperbolic metric) from the geodesic with same end
points on $\partial D$;

\item
$k_g=1$: horocycles, which are Euclidean circles internally tangent to
$\partial D$;

\item
$k_g>1$: hyperbolic circles, which are also Euclidean circles lying in the
interior of $D$.
\end{itemize}

A relation of this paper by Hedlund to magnetic geodesic flows was
pointed out by Ginzburg \cite{G}.

Any compact hyperbolic oriented two-manifold $M$ is the quotient of $D$
with respect to an action of
some discrete subgroup $\Gamma \subset SU(1,1)/\pm 1$:
$M = D/\Gamma$.

Let us consider the magnetic geodesic flow on $M$ with $F=d\mu$ on
the energy level $E < \frac{1}{2}$.
On $D$, the universal covering of $M$,
the pullback of a trajectory $\gamma$ is a
hyperbolic circle with the (hyperbolic)
center at $x_\gamma$. It is clear that,
given a fixed energy level $E$, to any point $x \in D$ there corresponds
a unique hyperbolic circle with constant geodesic curvature
$k_g = \frac{1}{\sqrt{2E}}$ centered at $x$.
Moreover to any trajectory on $D$
there corresponds its hyperbolic center $x_\gamma$
which is preserved by the flow. Transformations from $SU(1,1)/\pm 1$ maps
a hyperbolic circle into a hyperbolic circle preserving
the geodesic curvature and mapping the center into the center.
Therefore to any real-analytic function $f: M \to \R$
there corresponds a real-analytic integral of motion
$$
I_f(x,p) = f(x_\gamma)
$$
where $x_\gamma$ is the hyperbolic center of some pullback of a trajectory
$\gamma$ passing through $(x,\dot{x}) \in TM$. Here we identify $x_\gamma$
its projection in $M = D/\Gamma$ and $p_i = g_{ik}\dot{x}^k$.
Thus we have an integral of motion
$I_f$ defined on the set $\{E < \frac{1}{2} \}$.
This integral $I_f$ generates a Hamiltonian flow commuting with
the magnetic geodesic flow. This proves the theorem.

{\sc Remark} 2. Take a couple of real-analytic functions $f,g: M
\to \R$ which are functionally independent almost everywhere and
construct two independent integrals of motion $I_f$ and $I_g$.
These integrals do not commute with respect to Poisson brackets.
However fixing their generic values we single out isolated closed
trajectories of the flow.

{\sc Remark} 3. The additional integral from the proof of Theorem
depends analytically on the parameter $E < E_{\mathrm{cr}}$.


\begin{thebibliography}{99}

\bibitem{Hedlund}
Hedlund, G.A.
Duke Math. J. {\bf 2} (1936), 530--542.

\bibitem{CIPP}
Contreras, G., Iturriaga, R., Paternain, G.P., and Paternain, M.
Geom. and Funct. Analysis {\bf 8} (1998), 788--809.

\bibitem{T}
Taimanov, I.A.
Math. USSR-Izv. {\bf 30} (1988), 403--409.

\bibitem{BT}
Bolsinov, A.V., and Taimanov, I.A. Inventiones Mathematicae {\bf
140} (2000), 639--650.

\bibitem{Novikov}
Novikov, S.P.
Russian Math. Surveys {\bf 37}:5 (1982), 1--56.

\bibitem{G}
Ginzburg V.L.
Math. Z. 1996. V. 223. P. 397--409.



\end{thebibliography}
\end{document}